\documentclass[a4paper,12pt,leqno]{article}

\usepackage{geometry}
\usepackage{tikz-cd}\usepackage{bussproofs}
\usepackage{amssymb,amsmath}
\usepackage{proof}
\usepackage{psvectorian}
\usepackage{CJKutf8}
\usepackage[utf8]{inputenc}
\usepackage{graphicx}
\usepackage{graphicx}
\graphicspath{ {images/} }
\usepackage{qtree}
\usepackage{mathtools}
\usepackage{hyperref}
\usepackage{listings}
\usepackage{indentfirst}
\usepackage{titlesec}
\numberwithin{equation}{section}

\newtheorem{defn}[equation]{Definition}

\newtheorem{rem}[equation]{Remark}
\newtheorem{exm}[equation]{Example}

\newtheorem{notat}[equation]{Notation}
\newtheorem{newpar}[equation]{}

\newtheorem{xdefn}{Definition.}
\newtheorem{xproposition}{Proposition.}
\newtheorem{xcorollary}{Corollary.}
\newtheorem{xrem}{Remark.}
\newtheorem{xexm}{Example.}
\newtheorem{xlemma}{Lemma.}
\newtheorem{xtheorem}{Theorem.}
\newtheorem{xnotat}{Notation.}
\newtheorem{xnewpar}{\it}
\newtheorem{xproof}{{\it Proof. }}
\newtheorem{xproofof}{{\it Proof}}

\newenvironment{newparagraph*}[1]{\begin{xnewpar}\hspace*{-1.5mm}{#1}. \rm}{\end{xnewpar}}

\newenvironment{definition*}{\begin{xdefn}\em}{\end{xdefn}}
\newenvironment{remark*}{\begin{xrem}\em}{\end{xrem}}
\newenvironment{example*}{\begin{xexm}\em}{\end{xexm}}
\newenvironment{notation*}{\begin{xnotat}\em}{\end{xnotat}}
\newenvironment{proposition*}{\begin{xproposition}}{\end{xproposition}}
\newenvironment{corollary*}{\begin{xcorollary}}{\end{xcorollary}}
\newenvironment{lemma*}{\begin{xlemma}}{\end{xlemma}}
\newenvironment{theorem*}{\begin{xtheorem}}{\end{xtheorem}}

\titleformat*{\section}{\large\bfseries}

\begin{document}

\title{Modern incarnations of the Aristotelian concepts of Continuum and Topos}
\author{Clarence Protin\footnote{This work is financed by national funds through FCT – Fundação para a Ciência e a Tecnologia, I.P., within the scope of the project REF with the identifier DOI \url{10.54499/UIDB/00310/2020}. },  Centre of Philosophy of the University of Lisbon.}

\maketitle
\begin{abstract}
The aim of this paper is i) to argue for the feasibility and fruitfulness of a balance between the phenomenological method seeking intuitive evidence and the axiomatic-deductive method and ii) that there should be a mutual understanding between philosophy and mathematics and a cultivation of a historical self-awareness with regards to their common source in Greek philosophy.
To this end we show how Aristotle's theory of \emph{sunekhês, apeiron} and \emph{topos} and related notions can be given a rigorous interpretation in terms of modern topology and geometry as well as category theory. This is facilitated by the fact that in Aristotle himself we already find a balance between intuition and formal logic.  We also show how these powerful Aristotelian intuitions and concepts are found incarnated in diverse domains of modern mathematics.

\end{abstract}


\section{Introduction}

As Kant points out in his \emph{Critique of Pure Reason}, there are two main sources of human knowledge: intuition and concepts. In modern terms, we can seek to found knowledge by means of ultimate evidences (Husserl's \emph{Wesenschau})	or according to the axiomatic-deductive ideal pioneered by Leibniz and Frege. Although sometimes presented as conflicting tendencies (i.e. the intuitionist or constructivist vs. formalist strife in the foundations of mathematics) there is nothing inherently incompatible in these two approaches. Even Husserl allowed  an important role for an axiomatic-deductive treatment of certain aspects of his project of a Pure Logic (cf. his theory of 'manifolds'). An example of the  complementarity between these two modes of knowledge is  what pertains to space and time (or more generally the continuum) and their associated topological notions. 
In this paper we investigate some aspects of Aristotle's treatment of the continuum, its allied notion of infinity (\emph{apeiron}) as well as his theory of the \emph{topos}. We show that in Aristotle's theory there is a complementarity  between intuition and the axiomatic-deductive method and that a better understanding of Aristotle's theory is obtained (following the pioneering work of René Thom) by employing notions of modern topology and category theory. 
We also take a tour of some recent developments in topology and category theory and show how these developments can be seen as a return to, or rediscovery, of Aristotelian theory. More specifically, we argue that these developments represent an elaborate development of the intuitions and concepts present in an embryonic state in the Aristotelian theory. This elaboration is particularly strong on the axiomatic-deductive side. The desirability of a corresponding development from the intuitive phenomenological angle is patent, together with an awareness of an historical and philosophical continuity with the Aristotelian tradition. We are lead to conclude that philosophy and mathematics should be cognizant of each other and  not loose sight of their connection to a common historical origin.

This paper is organized as follows.
In section 2 we present the outline of Aristotle's theory of \emph{sunekhês, apeiron} and \emph{topos} as well as related notions. In section 3 we discuss how the concept  \emph{sunekhês} might be held to  play a powerful role in Aristotle's philosophy beyond a mere characterization of quantity.  In section 4 we give formulations of the key Aristotelian notions in terms of classical topology and sheaf theory. In section 5 we  review some examples of how Aristotelian notions incarnate in modern mathematics. In section 6 we lay the groundwork for a more adequate formalization of \emph{sunekhês} in terms of locales and topos theory and point out the significance of the former for Aristotle in general.  The final section ends with some conclusions about philosophical methodology, the relationship between philosophy and mathematics and the importance of historical awareness. An appendix is included in which are gathered some basic definitions from topology and category theory  used throughout the text.

\section{Sunekhês, apeiron and topos}

Aristotle's \emph{Physics} deals with the fundamental principles of nature (\emph{phusis}). Among these are change, magnitude\footnote{We will later discuss the significance of Aristotle's  distinction in (209b8) between
	\emph{diasteme} (extension) and \emph{megethos} (magnitude).}, place and time. Aristotle  argues that these are all \emph{continuous} quantities or continua. We will later argue, following René Thom,  that the importance of continuity for Aristotle exceeds  the quantitative aspects of physical being and extends to other ontological categories as well.

What does Aristotle mean by a continuum or \emph{sunekhês} ? 
His analysis of this concept relies on the allied  concepts  of point, limit (\emph{peras} which can also be translated by boundary or extremity), infinity (\emph{apeiron}), division, contact and \emph{contiguity}.

A striking feature of sunekhês is that it does not consist of points (it is not an aggregate or set of points) nor can it result from the joining together of points (231a21-b10)\footnote{Unless stated otherwise, all references are to the \emph{Physics} and we follow the text of \cite{Ross}}. This is surely one of the strongest known examples of a \emph{purely intuitively motivated theory}. In fact all the related concepts are equally intuitively evinced and form a tightly connected whole.

For Aristotle $\emph{peras}$ is intuitively the limit, boundary or extremity of an object.  Thus if a point has as limit it can only be the point itself. This also follows from the Euclidean definition of a point - that which has no parts - provided of course  we define 'part' as 'proper part'.

The passage (210a15-25) discussing the meaning of something being in something else for the purpose of establishing the definition of \emph{topos} can  be read as a proto-mereological discussion concerned with the general notion of parthood spanning different ontological categories. Reflexivity of the parthood notion is explicitly excluded, thus pointing to a restriction to proper parthood.

In (227a6-16) we get the following characterization of \emph{sunekhês}: the limits of its parts are united. Being united is a stronger notion than mere contiguity which depends in turn on the notion of place, \emph{topos}. Contiguous objects have limits 'in the same place'. This corresponds to 'contact'.

Since we cannot distinguish between the limits of a point and the point itself, the coming together (union) of two points will result in complete coalescence and
thus a continuum could never be constituted from points. 

Aristotle gives us in fact two fundamental complementary characterizations of \emph{sunekhês}:
(i) that whose parts have united limits (and hence cannot consists of points) and (ii) that which is \emph{infinitely}
divisible.  The core concept of infinity \emph{apeiron} is intimately bound up with that of $\emph{sunekhês}$. The characterization (i) may be termed 'connectivity', a Latin-derived term to which it is etymologically cognate.
Aristotle conceives a kind of correspondence or dependence between the different continua involved in change. For instance change itself, magnitude and time. This correspondence is expressed as: a division of one continuum must correspond to a division of the others. This principle is widespread in the proofs in Book VI of the \emph{Physics}. One of the most interesting proofs (237b23-238a20) is that there is no finite movement in an infinite time (i.e. unbounded infinite time). We will return to this proof later.

Aristotle distinguishes between an infinity which 'goes beyond any finite magnitude' (206b20) and a bounded infinity defined inductively by
successively adding  segments of half the previous length (this expresses the convergence of the geometric series).
In  illuminating passages (206b3-12, 207a7-25, 209b8) we see that such a bounded infinity is that which has a limit, 
not in itself but in something external to itself. It is something incomplete or incompleted, \emph{aoristos}(209b9),
imperfect or not yet perfect or whole (207a7-15). 
It can  be an expression both of a completable or incompletable process. It can express potentiality \emph{dunamis} and generated being \emph{oude menei(...)alla gignetai}, it does not rest(...)but  is in the process of becoming, and is likened to time (207b15).

There are in fact three basic kinds of \emph{apeiron} discussed in the \emph{Physics}. The boundless infinity exemplified by adding segments of the same length to each other and two kinds of bounded infinity exemplified as follows: one involving adding segments to each other in a fixed smaller proportion (as in the convergent geometric series) and the other (which we  call a 'vanishing infinity' ) obtained by successively dividing a given segment according to a fixed ratio. In Book VI Aristotle deploys such a vanishing infinity in proving that there is no first instant in which a change first occurred (236a7-26). All these infinities do not exist in act, but potentially, they exist as a processes, in particular as a discrete temporal processes.

An interesting aspect is that a bounded infinity can be used to define an object. Bounded infinities have, as we saw, a limit, not in themselves but exterior to themselves.  \emph{Thus a limit can be defined by a infinite process}.  Let us consider Aristotle's definition of \emph{topos} in this light. It is: \emph{to peras tou periekhontos somatos kath'ho sunaptei to periekhomeno}(212a5-6),
\emph{the limit of the enveloping body according to where it touches what it envelopes}.  In the final refinement of the definition the envelop is further required to be the first immobile  envelop (cf. the example of a boat in a river, where the river itself is the \emph{topos} of the boat). This strongly suggests that we can consider a bounded infinity consisting of all successively smaller enveloping bodies containing a given body. The limit of such a bounded infinity (restricted to immobile envelops) will be precisely the \emph{topos} of the body.

A body and its topos may or not be contiguous. Nor does the definition of topos employ the concept of contiguity. Thus there is no circularity in Aristotle's definition of contiguity as having limits in the same place (which corresponds to the fundamental notion of 'contact', the closest two objects can come together without loosing the individuality of theirs limits).

A careful reading of Aristotle's \emph{Physics} shows that \emph{sunehês, apeiron} and \emph{topos} and their allied concepts are not only intuitive and phenomenologically evident but the definitions given and properties proven are all conformity with (or at least the prototype of) the ideal of axiomatic-deductive clarity and rigour.

\section{Aristotle's philosophy of continuity}

Before moving on to discussing how the previously notions can be interpreted in terms of modern topology and category theory we find it advantageous to first go through briefly how continuity and its circle of intuitive notions play a larger role in Aristotle's philosophy than just as fundamental aspects of the category of quantity. We refer the  reader to seminal works of René Thom\cite{Esq} for further discussion and justification of this position. 

There can be little doubt that Aristotle was aware that the category of quality encompasses genera which can be continuously parametrized. Within these 'continuous' qualities some admit opposites (for instance wet-dry and hot-cold) and some do not.  A single body's state can be decomposed into a sequence of independent qualities, some admitting opposites, some not. In modern terms some qualities are represented by \emph{scalar} quantities (positive real numbers) and some by vector quantities (the real numbers  including negative reals - which in modern terms correspond to 'charge', or more generally vector spaces)\footnote{in Aristotle's proof that the world cannot be spatially infinite a postulate is employed that can be read as requiring the total canceling out of positive and negative charges.}.  Thus the four elements are characterized phenomenologically by four different regions in the two-dimensional vector space consisting of the product of the wet-dry and hot-cold genera each represented by a real line. Not only does the role of continuity transcend the domain of strict quantity but it has been suggested that much of modern physics was already implicitly present in Aristotle (cf. Thom's concluding remark in \cite{Bio}). 

A bolder and more controversial  thesis of Thom relates continuity to  \emph{hupokeimenon} (substance, which does not belong to the category of quantity but in which, according to Aristotle, quantity subsists)  as well to the whole intelligible realm of genera and species.

The continuity of the \emph{hupokeimenon} might be given the following justification via the consideration of the logical significance of limit (\emph{peras}) as a kind of interface, that which allows an entity
to relate to other entities or to the surrounding environment or be capable of receiving diverse predicates whilst remaining itself. These considerations have roots in Plato's \emph{Parmenides}. In the first
hypothesis of the  'dialectical' section,  the One is considered in its absolute simplicity.
It could not manifest (or exist)  for then it would have to be in a certain place.
But being in a place it would be \emph{enveloped} by that place, that is, its extremities would touch, be in contact, with the surrounding environment. Plato then continues to argue (in \emph{Parm}. 138A-B) that since  extremities are  parts of
the entity and the One is by definition an entity without parts,  the One
could not be in a place and hence could not 'exist'. According to Scolnikov \cite{Scolnikov}[pp.95 - 100] this problem is resolved in the second hypothesis
wherein the Existing One is considered as a complex unity, a composite (\emph{sunolon})
which allows participation (or predication) whilst preserving its underlying unity.
Scolnikov \cite{Scolnikov}[p100] makes the interesting
observation that this complex is in concrete instances always finitely determined or divided (\emph{diareton}), although
it can potentially be divided infinitely - \emph{diareisis} presupposes the possibility of making distinctions within a
form and requires, in principle, that division can go on indefinitely (infinitely).
But potentially endless divisibility is a chief property of the Aristotelian \emph{sunekhês}:(200b19-20) \emph{to eis apeiron diaireton sunekhês on} - the continuum is what is infinitely divisible.

Natural beings consist of both matter and form,
and for Aristotle form seems much like a continuum which is 'stratified' according to different locally homogeneous qualities.  This is at the heart of Aristotelian biological division
between \emph{homeomeres} and \emph{anomeomeres}.

As Thom argues in \cite{Esq}, scattered throughout the \emph{Physics} and \emph{Metaphysics} are passages suggesting Aristotle made a strong analogy (more, perhaps, than a significant metaphor) between the relationship of form and matter and the relationship of limit and (bounded) infinity. This was also extended to the analogy with the relationship between act and potentiality and even the relationship between a genus and its differences giving rise to species. As stated in the \emph{Metaphyiscs}: \emph{genos hos hulê},  genus is like matter.

\section{Modern mathematical interpretation of Aristotle' theory}


 We can either use modern mathematical concepts as a hermeneutic tool for Aristotelian philosophy or we can use aspects of Aristotelian philosophy to interpret certain concepts and results of modern mathematics (independently of whether the found correspondences are historically necessitated or represented independent convergence or rediscovery after an explicit rupture).
 
 There is also the problem of what we precisely mean by 'modern'. In our view a convergence with Aristotle if already present in classical 'modern' mathematics of the turn of the 20th century, only reaches its full force and significance with more contemporary developments, specially in category theory.
 
 The purpose of this section is to show how modern topology and category theory can be used as a hermeneutical tool to study Aristotle's theory of \emph{sunekhês, apeiron, topos} and related notions.
 
 In Aristotle's \emph{Physics} we find a balance and complementarity between intuitionism and axiomatic-deductive methodology\footnote{Thom himself as was only interested in Aristotle the proto-phenomenologist, dismissing completely the formal axiomatic-deductive Aristotle \cite{Esq}.  }. Modern mathematics on the other hand is characterized by a strong  focus on axiomatic-deductive methods (often involving an explicit rejection of anything epistemologically resembling Husserl's phenomenological method attributing an important role to intuitive evidence). Opposing schools took an equally strong view with regards to the rejection of the epistemic role of formal logic and axiomatic-deductive methods. This is even true of Husserl (specially in his later phase)  but not excluding an important and legitimate (albeit restricted) role for such methods. Be that as it may there can be no doubt that intuition plays a powerful role in the thought-processes of the ordinary 'working mathematician'. And it is both the formal definitions and intuitive content of some aspects of modern mathematics that we wish to deploy as a hermeneutical tool in this section.
 
 Although the idea of space not being constituted by points is an old one, as are attempts to elaborate a formal axiomatic-deductive theory of such theories, it has never been traditionally the framework of mainstream mathematics and theoretical physics. Its fullest and most sophisticated development took place within category theory, in particular through the work of Alexander Grothendieck and William Lawvere. But in order to motivate our use and discussion of some advanced category theoretic notions in relation to Aristotle, we find it convenient to start with basic notions from classical general topology\footnote{For a good introduction to general topology see for instance \cite{Kelley}}, which is still based on the idea of space as a collection of points.
 
 From here on we assume the reader to be familiar with the definition of a topology on a set (in particular the standard topology on the real line or Euclidean plane) and with the definition of open set, closed set, interior, closure, boundary, connectivity and continuous function.
 
 At first glance the two characterizations i) and ii) in Aristotle's definition of \emph{sunekhês} applied to a topological space seems to correspond to: (i) connectivity and (ii) to a certain decomposability property\footnote{this characterization is clearly related to a Hausdorff condition. In the one dimensional case this means that we can potentially divide any line, in particular the line connecting two points and so obtain
 	disjoint segments one containing each point.}.  The example of the coarse (or chaotic) topology shows that the two conditions are arguably independent.
 
 Before going into detail let us examine \emph{apeiron}.  From the analysis in the last sections we see that \emph{apeiron} is a complex and multi-faceted notion which spans several modern topological concepts. We posit the following two complementary characterizations of \emph{apeiron}:\\

 1)  Any set $A$ which does not contain at least one point in its closure $\overline{A}$ (for instance an open or semi-open interval).
 
 2)  A net $n$  which is either strictly increasing or strictly decreasing, that is either  for all $x\leq y$ with $x\neq y$, for $x,y$ belonging to a directed set $D$,  implies $n(x)\subsetneq n(y)$ or for all $x\leq y$ with $x\neq y$ implies $n(y)\subsetneq n(x)$.\\
 
 In 2) we do not require that $n(x)$ be an open set.  We observe that an \emph{apeiron} of the first type can be obtained by taking the union $\bigcup_{d \in D} n(d)$ (a particular case of limit) of a suitable net $n$. For instance Aristotle's construction corresponding to the convergent geometric series yields a semi-open interval. A point can likewise be seen as the limit of a net of strictly decreasing sets all containing that point. We will call such an \emph{apeiron} a \emph{vanishing infinity}\footnote{It is interesting to observe that these two kinds of limits corresponding to the two kinds of nets in the definition are reflected in the two fundamental kinds of Borel sets $F_{\sigma}$ and $G_{\delta}$, those that are countable unions of closed sets and countable intersections of open sets (the semi-open interval $(a,b]$ being an example of both). The same duality is furthermore present in the standard property of Borel measures wherein the measure of a set $E$ can be given equivalently as the infimum of the measures of open sets containing the set or the supremum of the measures of compact sets contained in the set.}

 As mentioned previously, Aristotles makes a distinction in (209b8) between
 \emph{diasteme} (extension) and \emph{megethos} (magnitude).  In modern terms the \emph{diasteme} would be the \emph{interior} (hence an open subset) of the megethos
 considered as a closed interval (cf. \cite{Pell}[p.207 Note 5] ).

 
  It follows naturally that the concept of \emph{peras} corresponds to the boundary of set as in 1) (for instance, the boundary of an open set)  or the limit (in the modern sense) of a net as in 2). But of course some restrictions would have to be made in the general setting to guarantee that we are working with intuitively justifiable  'tame' sets  in 1) and 2) so as to rule out such 'pathological' objects as the Cantor set.
  
  In the ordinary topology on the real line or Euclidean space a general open set can be quite a complex object. But the fact of the matter is that all such open sets can be seen as \emph{generated} from (i.e. in the sense of a basis of a topology) more basic elementary open sets, for instance open intervals or balls. Such basic open sets have very Aristotelian characteristics. They function like local homogeneous patches that generate (by infinite unions) all the open sets of the topological space. They express the concept of the local spatial quality around a given point in which size  is irrelevant (cf. Hegel's concept of indifferent quantity in \cite{Hegel}). Metaphorically a (basic) open set corresponds to matter (\emph{hule}) and potentiality (\emph{dunamis}). 
  
  Giving a basis for a topology is the same as giving at each point a vanishing infinity (i.e. a net whose limit is that point). The homogeneous nature (any local open part has the same quality as the whole) of the continuum is due to it being generated by the same archetypal vanishing infinity. Also a topology cannot be generated by a finite net nor identified with any fixed vanishing infinities for each point, because many different nets whose limit is that point give rise to the same topology. A basis around a point is a mediator which allows a point to be 'in' a continuum (something which it cannot do in itself).

  A closed interval or magnitude \emph{megethos} can be decomposed into  its interior \emph{diasteme} and its boundary \emph{peras}.  A point is a closed set and any point inside a magnitude divides it potentially into two separate magnitudes, it being the common limit of the interiors of these two magnitudes.
  On the other hand two magnitudes become a single magnitude by fusion and identification of limit-points. Thus the closed point as a limit mediates both division and junction.
    More generally the genesis of an entity can be seen as a separation
 from a matter by means of a cut along its boundary, thus  defining and giving form to the entity. 
 
 In the treatment of motion in Book VI there are many (semi-)open sets or in general \emph{apeiron} objects  defined: for instance the interval of time in which change takes place is open while the magnitude of completed change is always closed (and bounded).
 
 In Euclidean space sets homeomorphic to three-dimensional balls express the extension of real, actually existing beings. Whilst for instance surfaces have a dependent, incomplete existence.  In general this suggests that substances correspond to closed sets with non-empty interior (and bounded, thus compact, if we accept with Aristotle the finitude of the universe). Thus this is the example of the most ontologically robust \emph{sunekhês} alongside closed intervals of time or closed intervals of spatial locomotion. We may wish to indentify this with \emph{megethos}; however Aristotle apparently does consider an unbounded \emph{megethos} in Book VI (the time starting from a given moment).  The interior of the above kind of \emph{subekhês} expresses the matter and is a type 1) \emph{apeiron}. The boundary is its limit and its form.
 
 The important task that remains is formulating conditions i) and ii) in the setting of general topology (we will in fact do this later for a more general concept of topology without points which is more faithful to Aristotle's own theory).

 Compactness is significant in relationship to Aristotle's proof
 that there can be no infinite motion in a finite time (237b23-238a11-20) or finite motion in an infinite time. We mentioned how Aristotle lacked the modern concept of function; rather he worked with two-way correspondences between various quantities such as motion, magnitude and time.  These were special kinds of correspondence in that divisions of one quantity corresponded to divisions of the other. Thus if we divide an interval of time we get a division of the corresponding motion and corresponding magnitude traversed, the parts corresponding respectively. Also if we divide a motion we get corresponding divisions in the magnitudes traversed and times taken, the parts corresponding respectively. It would be  interesting to investigate how this relates to the modern definition of a \emph{continuous} function.  Perhaps such correspondences can be captured by homeomorphisms or pairs of continuous functions (or, as discussed subsequently, in a general categorical setting as adjunction pairs on sites). Anyhow, if we consider, as is natural, that such correspondences are mediated by continuous functions then Aristotle's conclusions are valid. The image of a compact set must be compact. Thus there can be no infinite motion in a finite time nor finite motion in an infinite time (because the image via a homeomorphism of $[0,\infty)$ cannot be a closed interval). Aristotle's proofs use a kind of 'covering' which  recalls the modern definition of compactness in terms of open coverings.
 
 We point out that the modern definition of continuous function can in certain circumstances be given a definition in terms of the more \emph{sunekhês}-like property of connectivity. The image under a continuous function of connected set is connected. A particular case of this property is Bolzano's intermediate value theorem (taught in elementary calculus courses)
 which expresses \emph{natura non datur saltus},  that a continuous function has no jumps or holes,  cannot skip intermediate values (or that the graph is a connected curve in the plane). Indeed the continuity of a function can be defined more intuitively in terms of the arcwise connectedness of its graph.
 
 But for Aristotle the world is not completely governed by what Leibniz and we today would call continuous processes. A passage (228a26-30) attests how Aristotle recognizes that globally change is not
   continuous, only locally continuous, a concatenation of \emph{contiguous} continuous changes.
 
 Within the framework of classical point-set topology let us see how the concept of \emph{topos} could be formulated. For Aristotle the notion of topos involves real qualified entities in nature, not geometric forms. As Léon Robin 
 wrote in\cite{Robin}[p.142]: a topos is 'always qualified'. It is not hard to see that if we define the topos of a 'nice' set $B$ in space to be the limit of the net of (basic) open sets\footnote{or closed or more general sets.} $U$ containing it (that is, the boundary of the intersection of all such open sets) then we get simply the boundary of $B$ itself, since the intersection of all such sets is just the closure of $B$. But this contradicts the necessary separability of container and contained, the difference between the form and topos of an entity.
 
 But if we consider the set $B$ endowed with a continuous field of qualities then Aristotle's definition can be given a precise definition in modern terms. We assume that there is a space of qualities $\Phi$ and that the entity corresponding to a closed set $B$ lies within a larger environment given by, let us say, Euclidean space $E$ with the usual topology. For every open set $O$ of $E$ there are functions $s: O \rightarrow \Phi$, the set of such functions being denoted by $\Phi(O)$. These represent the possible qualities of the region $O$. This is a particular instance of the concept of sheaf on a topological space \footnote{see the appendix for the definition of the modern mathematical concept of sheaf. An introduction  to Sheaf Theory can be found in \cite{Moerdijk, Godement}.}. The question is now: what are the qualities on $B$ and how do they relate to the sheaf of possible qualities on $E$ ? The sheaf is only defined on open sets of $E$ not on $B$ which is not open. We can give $B$ the induced topology but it is still not immediately clear how to define a $s: V \rightarrow \Phi$ on one of the open sets $V \subset B$ of the induced topology. In particular how can we define the set of possible phenomenological functions $s: B \rightarrow \Phi$ on $B$ itself ? There is no other way than to consider $\Phi(U)$ for each open set $U$ containing $B$. We consider the \emph{apeiron} of the net of nested open sets $U$ and their qualities $\Phi(U)$ and take the \emph{limit} (in the modern sheaf-theoretic sense) of the $\Phi(U)$.  This is well-defined and captures perfectly the idea of the space of qualities of $B$ relative to its \emph{immediate} environment in which the relationship to the environment is mediated precisely by the topos of $B$ in $E$ specified to the sheaf $\Phi$. Thus we can give at last an interpretation to (211b10) 
 we have \emph{en tauto gar ta eskhata tou periekhontos kai tou periekhomenou}, 
 \emph{the limits of the container and contained 
 	coincide} (211b12-14), but without fusing and being identical.

  While we saw that while just considering \emph{apeiron} of the open sets containing $U$ is vacuous, once we consider a sheaf of qualities it becomes highly significant and non-trivial.  A similar construction can be carried out for any $V$ of the induced topology of $B$. In sheaf theory this corresponds to a general construction called the \emph{inverse image sheaf} which is well behaved for sufficiently 'tame' kinds of set $B$. This construction is also a generalization of the infinitesimal calculus. The derivative of a function at given point cannot be calculated from the value of the function at that point alone. A certain neighbourhood of values is required.  In fact if we consider $B$ to be just a point $p$ in $E$ then the previous construction yields what is called the \emph{stalk} of $\Phi$ at $p$, the qualities around an indefinitely small neighbourhood of $p$. The Aristotean concept of (vanishing) \emph{apeiron} in this interpretation justifies the 'flowing' (cf. Newton's \emph{fluxiones}) nature of the infinitesimal.  Of course if we just consider open sets containing $p$ then the limit is $p$ itself. Cf. (209a7-13) where the point is stated to be identifical to its topos. Waterfield \cite{Water}[p.253] also comments on a connection to the first hypothesis of the Parmenides.

   \section{Some mathematics reflecting Aristotelian intuitions and concepts}
   
   In this section we take the opposite approach. We examine some modern mathematical concepts and point out how they can be given an Aristotelian interpretation or clearly reflect Aristotelian concepts even when it is almost certain that there was no direct or conscious influence.
   
   The Aristotelian topos is found in algebraic geometry in the guise of the formal completion of a subvariety $Y$ in a variety $X$.
   In \cite{Hart}[p. 190] Hartshorne writes:\\
   
   \emph{The formal completion of $Y$ in $X$ (...) is  an object which carries information about all the infinitesimal neighbourhoods $Y_n$ of $Y$ at once. Thus
   	it is thicker than any $Y_n$, but is  contained inside any actual open neighborhood of $Y$ in $X$.}\\

   Elementary algebraic topology also offers us many striking embodiments of the Aristotelic topos, situations in which we are interested in the behaviour in the neighbourhood of a boundary. This is the main intuition behind the concept of relative homology of a pair $(X,A)$ and the excision theorem. According to \cite{Vick}[p.44] : 'In particular a chain in X is a cycle modulo A if its boundary is contained in A. This reflects the structure of $X-A$ and the way that it is attached to A. In a sense, changes in the interior of A, away from its boundary with $X-A$, should not alter these homology groups'.
   		A very useful property of CW-complexes is that they are locally strong deformation retracts around the 'place' in which the cells are attached.

   		 In continuum mechanics \emph{Cauchy's stress principle} is related to the 'act' of an infinitesimal homeomeric volume
   		 of the substance on its immediate environment through the defining region of the act, the surface: the \emph{surface force}.
   		 The infinitesimal volume will always be dynamic and have a surface; Cauchy's principle states 
   		 that we can take the limit of the force/area ratio of the force acting on a given small area of the surface determined by a
   		 normal direction because the moments
   		 will vanish at the limit. Considering each direction we obtain a (symmetric) tensor, the stress tensor, at each infinitesimal point.
   		 Cauchy's principle is an example of the dynamical  nature of the infinitesimal.

   		One of René Thom's most important papers in topology\cite{Strat} concerns stratified sets and morphisms. The theory involved - which deals with a generalisation of smooth manifolds and algebraic and semi-algebraic spaces -  dates back to the work of Hassler Whitney (1907-1989) on singularity theory .   It can also be seen as a development of the topological intuition present in the Aristotelic concept of topos. The topos is found in the way in which in a stratification a given stratum $X$ is related to  the adjacent strata $Y$ which are contained in its closure $\bar{X}$. This induces a decomposition of the boundary of $X$ into such adjacent strata which can be seen as forming the \emph{topos} of $X$. These strata $Y$ in general have a very complex topological relationship to $X$. Whitney's conditions (A) and (B) are imposed in order to maintain a tame, intuitive, non-pathological situation,  ruling out what can happen with spirals or  the crossings of the Whitney umbrella. These Whitney stratifications came closer to the intuitive notions of topos due to the smooth way in which the adjacent strata envelop and contain $X$. Many more mathematical examples that could be adduced regarding  the importance of the behaviour on the 'boundary' of object such as for example the Dirichlet problem  in partial differential equations in which a function on a  given region is completely determined by its values on the boundary of that region. 
   		
   		 The concept of a \emph{constructible
   			sheaf}\cite{Kashiwara}[p.320] involves a stratification of the underlyng space into locally homogenous qualities in a way
   		the captures admirably the Aristotelian concept of \emph{anomeomere}. Roughly speaking a locally constant sheaf describes qualities which are locally the same (homogenenous) but globally may conceal non-trivial information. For instance the earth appears flat locally but globally is endowed with spherical curvature. A constructible sheaf is a sheaf where its underlying topological space admits a stratification such that the restriction of the sheaf to each stratum is a locally constant sheaf. For instance the different kinds of cells or tissues in a living organism are organized according to geometric shapes  each having similar local genetic properties (they are homeomeres). See \cite{Esq} for a detailed discussion.

      Finally the theory of molecular toposes discussed in the next section is a striking example of a manifestation of \emph{sunekhês}: a theory of space based on open sets rather than points and whose properties are like Aristotle's characterization of \emph{sunekhês} defined in terms of connectivity and decomposition.

   \section{Category theory and sunekhês}

. 
  
  The Aristotelian theory of \emph{sunekhês} posits that continuous quantity is not an aggregate or set of points.   Thus the best modern mathematical interpretation, unlike the one given in the previous sections, must be built on similar suppositions. As we mentioned before, it is within category theory that such a theory has been carried out with the most sophistication as well as relevance to other areas of logic, mathematics and science. There are also independent philosophical reasons  why such a theory of continuous quantity should be preferred. Open sets are given to us directly and immediately via intuition, while points are constructed. And if we posit points or sets of points are being fundamental, then open sets, any topology in fact, is something arbitrary and constructed, something we impose on reality. 
  
  In this section we focus on the theory of locales and topos theory. These approaches are relational. In locale theory open sets are now themselves primitive entities rather than the elements of the set $X$ on which the usual notion of topology is defined. What counts now is the relationship between these open set primitives. In locale theory 'point' is a derived notion. There are locales which are non-trivial but hardly have any 'points'. There are others which are said to 'have enough points' which can be proven to be equivalent to the class of so-called 'sober' topological spaces. Thus locale theory does not involve a necessary rejection of points. Also we can naturally conceive of a theory which posits both open sets and points as primitive notions.
    
   A locale is a complete distributive
  lattice (which hence has a top element $1$ and bottom element $0$) and is a way of defining topology without points (see chapter II of \cite{Pedicchio} for an introduction). We might attempt to
   capture a point-free description of the Aristotelian \emph{sunekhês} properties i) and ii) exemplified in the standard topology of Euclidean space by the following axioms:
  
  \[ \neg \exists u,v. (1 = u\vee v)\, \& \, (u\wedge v = 0) \text{          $\quad \quad$   Global Connectivity}\]
  \[ 1 = \bigvee \{ z: \neg \exists u,v. (z = u \vee v) \&\, (u\wedge v = 0)\} \quad \text{   Local Connectivity}    \]
  \[ \forall u \exists w,v. ((w\wedge v = 0)\, \& \, \forall z. Con(z) \&(w\vee v \leq z) \rightarrow  u\leq z) \text \quad  \text{Decomposability}\]

  where $Con(z)\equiv \neg \exists u,v. (z = u \vee v), \&\, (u\wedge v = 0)$. We call a locale satisfying the axioms above an \emph{Aristotelian locale}.
  Local connectivity means that we can cover any open set with connected open sets and Global connectivity means that the whole space is connected. These axioms were obtained from  considering Aristotle's properties i) and ii). Property i) of course does not mean that any two parts have to have common limits (for just take two disjoint subintervals) rather that if they are contiguous (or share the same topos) then they share their limit. In view of ii) this means that given a decomposition as guaranteed by the Decomposability axiom and in which $w$ and $v$ are connected it cannot be the case that in fact $z = w \vee v$. Interestingly this can be seen as a kind of 'completeness' condition anticipating Cauchy and Dedekind. 
  
  Let us move on to topos theory and see if we can use it to give a definition of \emph{sunekhês} based on the above considerations for locales. Sheaves are usually defined in terms of the (opposite) category of open sets of a given topology. But we can define sheaves directly on locales or even on general categories.  Grothendieck's philosophy\cite{SGAIV} is that the category of sheaves on a space gives a more fundamental description of a space than its topology. A Grothendieck topology is precisely a way of defining a 'topology' on an arbitrary category without the need of an underlying set of points.
  The key concept is  that of a 'cover' which can be traced to the ancient Greek notion of a quantity measuring another quantity and which plays a role in many proofs in Aristotle's \emph{Physics}.  We thus can define sheaves
  on categories endowed with a Grothendieck topology (called \emph{sites}). Sheaves on a site in turn form a category called a \emph{Grothendieck topos}.
  A still more general concept of Topos, called an Elementary Topos (characterized by a small and elegant set of axioms) was discovered independently by William Lawvere in his research into the foundations of mathematics and intuitionistic set theory. Toposes can be used as a  semantics for formal logics in which we wish to consider a more general set of truth-values $\Omega$ endowed with an algebraic structure. The case of the set $\{true, false\}$ with the standard logical operations is the simplest Boolean algebra. Grothendieck topos semantics allow us to a have temporal or spatially localized version of truth and to give meaning
  to  sentences  such as 'The cat is black and white' and to interpret a relativized version of the principle of non-contradiction which incidentally is the one stated in Aristotle's \emph{Metaphysics}. But we are not proposing that Aristotle's logic or metalogic are intuitionistic logic\footnote{this is certainly contradicted by the logic of the \emph{Analytics}. See for instance \cite{pro} for a formalization of Aristotle's syllogistic.} only that intuitionistic logic or metalogic can be considered as locally valid for certain regional ontologies pertaining to Aristotelian physics.
  
   The poset of subobjects $Sub(A)$ of an object $A$ in a category can have diverse properties depending on the type of category. In the case of a topos they are Heyting algebras. $Sub(A)$ is a kind of generalised space.
  
  We wish to define an Aristotelian topos as a topos satisfying the analogues of the axioms for Aristotelian locales.  In\cite{DiaconBarr} we are offered a definition of molecular (also called locally connected) topos.  This means that every object $E$ can be written as a categorical sum $E = \Sigma_i M_i$ where $M_i$ is a molecule (or connected object) 
  meaning that $M_i$  cannot be written as a categorical sum of two proper subobjects. 
  This is the precise analogue of our Local Connectivity Axiom. There is also an analogue of the Global Connectivity Axiom.
  Interestingly in \cite{LoCon} we find that we can characterize connected objects in a topos $T$ in terms of their associated hom-functor $hom_T(A,-)$
  preserving finite coproducts (this appears to have been first discovered by Grothendieck in \cite{SGAIV}). And it turns out
  that the topos of sheaves over a topological space $X$ is locally connected if and only if $X$ is locally connected in the classical sense.
  
  Thus it would be desirable to extend the to toposes the remaining Divisibility Axiom for Aristotelian locales in order to obtain a candidate for a topos theoretic formulation of \emph{sunekhês}.

  There is a relationship between intuitionistic logic and general topology which is suggestive for a topological reading of Aristotle. We stress again that this is not meant to imply that either Aristotle's logic or metalogic are  intuitionistic.
   Heyting algebras, the key algebraic component of the logical structure of a topos, have illuminating connections to classical topology and the concept of boundary.  If we are on the boundary of something, are we 'inside' or 'outside' that object ? A natural answer is : neither. Thus the 'outside' of an object \emph{must exclude the boundary}. This 'outside' is thus distinct from the set-theoretic complement which gives the set of subsets of a given set $X$  the structure of a Boolean algebra. But let us consider a topology $\tau$ on $X$ instead of all subsets. Then the complement $O^c$ of an open set is in general not an open set. The correct definition of 'outside'  of $O$ is that of the topological \emph{interior} of the complement $(O^c)^\circ$. In a similar way we can define operations which turn $\tau$ into a Heyting algebra in which $\vee$ corresponds to union and $\wedge$ to (finite) intersection, $\top$ to $X$ itself and $\bot$ to $\emptyset$. The failure of the law of the excluded middle in general, which reads in this case $O \cup (O^c)^\circ \neq X$, is seen as a direct consequence of the exclusion of the boundary from both $O$ and 
  $(O^c)^\circ$.
  
  In section 3 we discussed the central role of \emph{sunekhês} in Aristotle's general philosophy. Topos theory is an alternative to set theory which also allows us to construct alternative models of mathematics which satisfy or violate certain axioms such as the Axiom of Choice or the Continuum Hypothesis. To obtain a topos theoretic formulation of Aristotelian concepts it is relevant to try to find a topos in which, when considered as a semantics for its own internal language, all functions are continuous or at least locally continuous. 
   The construction of the reals via Dedekind cuts on the rational numbers can be generalised 
   to the topos of sheaves over a topological space $X$.
   The real number object in this topos corresponds surprisingly to 
   the sheaf of continuous functions over $X$. The objects are merely sheaves over $X$ (there is no continuity or real numbers involved in general). The real numbers are defined via a topos generalization of Dedekind cuts. And yet the sheaf of \emph{continuous} functions to the reals emerges. Somehow the generalized definition of real number object and the sheaf condition on a topological space combine to produce continuity. If this topos is considered as a model for mathematics it can be proven that all functions are continuous\cite{Moerdijk}[p.324], that is elements of $\mathbb{R}^\mathbb{R}$ satisfy the internal version of the continuity condition. It would be however more faithful to Aristotle if all functions were instead merely locally continuous and, of course, if the topos were not based on sheaves over a classical topological space. Thus we propose that this construction be studied for the general case of a Dedekind-type real number object  in a topos of sheaves over a locale or site.
   
   We saw that the analogue of 'continuous function' for Aristotle is a kind of decomposition-preserving correspondence between various quantities involved in motion: motion itself, magnitude traversed, time, etc.  This is also a prototype of a definition of continuous function that does not involve sets of points. In topos theory the concept of continuous function is  subsumed by that of \emph{geometric morphism}. One definition of a geometric morphism is: a functor between two toposes that preserves finite limits and has a right adjoint. Thus a geometric morphism is a pair of functors each going in the opposite direction (not unlike Aristotle's correspondences just mentioned). It is a purely topos theoretic concept which does not involve any topological space or set of points.  However a continuous map $f:X \rightarrow Y$ between Hausdorff topological spaces gives rise to a geometric morphism between the corresponding toposes of sheaves and all geometric morphisms arise that way\cite{Geomor}. Furthermore this situation can be generalized from maps of topological spaces to functors between sites which preserve covers. This seems a good approach to formalizing the decomposition-preserving correspondences of the \emph{Physics}.

  \section{Conclusion}
  
  We have seen that the Aristotelian concepts of \emph{sunekhês} and \emph{topos} and their circle of related notions express very powerful and deep-rooted intuitions. We showed that they are not only capable of being studied in a rigorous axiomatic-deductive way, using modern mathematics, but their power is attested in the way they emerge spontaneously in various domains of modern geometry and topology. 
  
  Our goal has been to illustrate the possible harmony and mutual enrichment between the phenomenologist's preoccupation with intuitive evidence (in this case geometric and topological intuition) and the formal rigour of the axiomatic-deductive method. We propose that these two methods already found a good balance in Aristotle.  It is certainly no coincidence that powerful and persistent Aristotelian topological intuitions also have received some of the most perfect and sophisticated axiomatic-deductive treatments to date.
  
  Aristotle's positing that the infinite in unknowable (for instance in his refutation of Anaxagoras in Book I of the \emph{Physics}) and at the same time only existing potentially is significant.  As a process and at an incomplete stage  the infinite is intuitively evident and knowable. So too is the limit (if it exists) of the infinite. It is only as a completed totality that it is unknowable. In classical set-theoretic extensionalist mathematics we reach in fact similar conclusions, for instance, that there is an uncountable set of indefinable real numbers which we cannot ever know individually or distinguish within such a set. Aristotle's position is rather one of epistemic optimism with regard to the progress of the power of reason grounded both in intuitive evidence and logic.
  
  Another closely related goal has been to show the fruitfulness that results from the shared knowledge of philosopher and mathematician and the importance of developing a philosophical and historical self-understanding  with regards to the richness of their  common source in classical Greek philosophy. Until the days of Husserl  it was normal for a philosopher to have some knowledge of mathematics and physics and to not only reflect upon the foundational questions of mathematics and physics but to see these questions as important for epistemology, ontology and the philosophy of mind as they are for the philosophy of nature.  The subsequent divorce of science and mathematics from philosophy has had the tragic consequence that the independent technical development of the latter
  has often led to it missing out on a mutually enriching interplay with the philosophical thought of the past.
  
  Thus we propose that further progress in both philosophy and mathematics will involve the simultaneous development of the phenomenological  (in particular as applied to geometric and topological intuition ) and the axiomatic-deductive methods as well as engagement in the interplay which spontaneously emerges between both.

  \section*{Appendix}
  
  We give here some definitions of a few important concepts in modern topology and category theory.
  
  Let $X$ be a set. Then a \emph{topology} over $X$ is a collection $\tau$ of subsets of $X$ such
  that $X$ and the empty set belong to $\tau$, the union of any collection of sets in $\tau$ belongs to $\tau$
  and the finite intersection of sets in $\tau$ is still in $\tau$.  The elements of $\tau$ are called
  \emph{open} sets. A \emph{closed} is the complement of an open set. A \emph{basis} for a topology $\tau$ is a collection $\mathcal{B}$
  of subsets of $X$ such that every open set can be expressed as a union of sets in $\mathcal{B}$.
  A basis for the standard topology of the real line is given by the open intervals $(a,b)$.  Given a point $x \in X$
  a \emph{neighbourhood} of $x$ is an open set containing $x$.
  Given a  subset $Y$ of $X$ the $\emph{interior}$ of $Y$ is the set of all points that have a neighbourhood contained in $Y$. The
  $\emph{closure}$ of $Y$ is the intersection of all open sets containing $Y$. It is always a closed set.
 The correspondence which associates to each subset $Y$ its closure is an example of a \emph{closure operator}.   For instance, the closure of the open interval $(a,b)$
  is the closed interval $[a,b]$. The \emph{boundary} of a set is the set of all points $x$ that satisfy the condition: any neighbourhood of
  $x$ will contain points both in $Y$ and not in $Y$. There is a general concept of \emph{compactness} but in
  the case of the ordinary topology on the real line (or Euclidean space $R^n$) it corresponds to being both closed and bounded (all points
  at a finite distance from the origin).
  
  A \emph{directed set} $D$ is a partially ordered set such that if $d_1, d_2 \in d$ then there is a $d_3$ such that $d_1 \leq d_3$ and $d_2 \leq d_3$. 
  A \emph{net} on a space $X$ is a map $n: D \rightarrow \mathcal{P}X$ for a directed set $D$. Here $\mathcal{P}X$ is the set of parts of $X$. The standard definition of net involves taking elements of $D$ into open sets of some topology $\tau$ on $X$. For the definition of \emph{limit} of a net see \cite{Kelley}.

  A \emph{category} can be thought of as a collection of \emph{objects} and for any two objects (possibly equal) a (possibly empty) collection of \emph{arrows} (or morphisms).
  Each arrow has a source and a target. If the target of an arrow $f$ is equal to the source of an arrow $g$ then we can compose
  them to form a new arrow $g \circ f$ which will have the source of $f$ and the target of $g$. Given an object $A$ there
  is a unit arrow $u_A$ which goes from $A$ back to $A$.  The formal definition of category requires an associative condition on
  composition and that the unit arrows function as units for composition, for instance $f \circ u_A = f$ when $f$ has source $A$.
  Given two categories $C$ and $D$ a \emph{functor} is a correspondence between objects of $C$ and objects of $D$ and arrows of $C$ and
  arrows of $D$ which preserves composition and unit arrows.
  Examples of categories are the category $Set$ of sets whose objects are sets and morphisms are maps between sets and
  the category $\mathcal{O}(X)$ of open sets of a topology on $X$ whose objects are the open sets for any two open sets $U$ and $V$ there is either a single morphism 
  if $U$ is contained and $V$ or no morphism at all if otherwise. In the category  $Set$ we can form the \emph{cartesian product} $X\times Y$ of
  two sets $X$ and $Y$ as well as the \emph{disjoint union} $X \sqcup Y$. The properties of these constructions can be
  generalised to obtain the concepts of product and coproduct for any category (and more generally the concepts of limit and colimit). A \emph{Topos} can be understood as
  a category which possesses generalizations of many fundamental properties of $Set$.
  
  Given a category $C$ we can form the \emph{opposite category} $C^{op}$ which is obtained by inverting the direction of all the arrows.
  A \emph{presheaf} is simply a functor $F$ from $C^{op}$ to $Set$.  Given an object $A$ of $C$ the elements of $F(A)$ are called \emph{sections}.
  If $U \subset V$ then this arrow corresponds to an arrow $F(V)\rightarrow F(U)$ (note that source and target are switched because we are considering
  the opposite category) called the \emph{restriction} map.
  Presheaves can be made into a category with a suitable definition of  morphism (natural transformations between functors).  If $C$ is the category of open sets of a topology on a space $X$
  then a \emph{sheaf} over $X$ is a presheaf on $C$ which satisfies the \emph{gluing} condition. This says that if we express
  an open set $U$ as a union (cover) of open sets $U_i$ and we consider sections $s_i$ on each $U_i$ such that they agree on all the intersections $U_i \cap U_j$
  then there is unique section on $U$ such that the sections $s_i$ are given by the restrictions corresponding to the inclusions $U_i \subset U$.
  An example of a sheaf is the correspondence which associates to each open set $U$ of a topological space $X$ the set of \emph{continuous
  functions} defined over $X$. An example of a presheaf which is not a sheaf is given by constant functions on open sets. Considering
  two disjoint open sets it is easy to see why the gluing condition fails.
  The concept of \emph{Grothendieck} topology is a generalization of the concepts involved in this property to an arbitrary category.
  The idea of covering $U$ by the  open sets $U_i$ becomes the concept of \emph{sieve}.

								\end{document}